\documentclass[titlepage,12pt]{article}

\usepackage{graphicx}
\usepackage{amssymb}
\usepackage{psfrag}

\newcommand{\R}{\mathbb R}

\begin{document}

\baselineskip=18pt

\begin{center}
{\Large{\bf Lyapunov Coefficients for Degenerate Hopf Bifurcations}}
\end{center}

\vspace{1cm}

\begin{center}
{\large Jorge Sotomayor}
\end{center}
\begin{center}
{\em Instituto de Matem\'atica e Estat\'{\i}stica, Universidade de
S\~ao Paulo\\ Rua do Mat\~ao 1010, Cidade Universit\'aria\\ CEP
05.508-090, S\~ao Paulo, SP, Brazil
\\}e--mail:sotp@ime.usp.br
\end{center}
\begin{center}
{\large Luis Fernando Mello}
\end{center}
\begin{center}
{\em Instituto de Ci\^encias Exatas, Universidade Federal de
Itajub\'a\\Avenida BPS 1303, Pinheirinho, CEP 37.500-903, Itajub\'a,
MG, Brazil
\\}e--mail:lfmelo@unifei.edu.br
\end{center}
\begin{center}
{\large Denis de Carvalho Braga}
\end{center}
\begin{center}
{\em Instituto de Sistemas El\'etricos e Energia, Universidade
Federal de Itajub\'a\\Avenida BPS 1303, Pinheirinho, CEP 37.500-903,
Itajub\'a, MG, Brazil
\\}e--mail:braga\_denis@yahoo.com.br
\end{center}

\vspace{0.5cm}

\begin{center}
{\bf Abstract}
\end{center}

\vspace{0.1cm}

In this paper are studied the codimensions one, two, three and four
Hopf bifurcations and the pertinent Lyapunov stability coefficients.
Algebraic expressions  obtained with computer assisted calculations
are displayed.

\vspace{0.1cm}

\noindent {\small {\bf Key-words}: Lyapunov coefficients, degenerate
Hopf bifurcation.}

\noindent {\small {\bf MSC}: 70K50, 70K20.}

\newpage

\section{Lyapunov coefficients}\label{S3}

\newtheorem{teo}{Theorem}[section]
\newtheorem{lema}[teo]{Lemma}
\newtheorem{prop}[teo]{Proposition}
\newtheorem{cor}[teo]{Corollary}
\newtheorem{remark}[teo]{Remark}
\newtheorem{example}[teo]{Example}

The beginning of this section is a review of the method found in
\cite{kuznet} and in \cite{kuznet2} for the calculation of the first
and second Lyapunov coefficients. The calculation of the third
Lyapunov coefficient can be found in \cite{smb2}. The calculation of
the fourth Lyapunov coefficient has not been found by the authors in
the current literature. The extensive calculations and the long
expressions for these coefficients have been obtained with the
software MATHEMATICA 5 \cite{math} and the main computational steps
have been posted in the site \cite{mello}.

Consider the differential equations
\begin{equation}
{\bf x}' = f ({\bf x}, {\bf \mu}), \label{diffequat}
\end{equation}
where ${\bf x} \in \R^n$ and ${\bf \mu} \in \R^m$ are respectively
vectors representing phase variables and control parameters. Assume
that $f$ is of class $C^{\infty}$ in $\R^n \times \R^m$. Suppose
(\ref{diffequat}) has an equilibrium point ${\bf x} = {\bf x_0}$ at
${\bf \mu} = {\bf \mu_0}$ and, denoting the variable ${\bf x}-{\bf
x_0}$ also by ${\bf x}$, write
\begin{equation}
F({\bf x}) = f ({\bf x}, {\bf \mu_0}) \label{Fhomo}
\end{equation}
as {\small
\begin{eqnarray}\label{taylorexp}
F({\bf x}) = A{\bf x} + \frac{1}{2} \: B({\bf x},{\bf x}) +
\frac{1}{6} \: C({\bf x}, {\bf x}, {\bf x}) + \: \frac{1}{24} \:
D({\bf x}, {\bf x}, {\bf x}, {\bf x}) + \frac{1}{120} \: E({\bf
x}, {\bf x}, {\bf x}, {\bf x}, {\bf x}) + {\nonumber} \\
\frac{1}{720} \: K({\bf x}, {\bf x}, {\bf x}, {\bf x}, {\bf x}, {\bf
x})+ \frac{1}{5040} \: L({\bf x}, {\bf x}, {\bf x}, {\bf x}, {\bf
x}, {\bf x}, {\bf x}) +  \frac{1}{40320} \: M({\bf x}, {\bf x}, {\bf
x}, {\bf x}, {\bf x}, {\bf x}, {\bf x},{\bf x})  \\+
\frac{1}{362880} \: N({\bf x}, {\bf x}, {\bf x}, {\bf x}, {\bf x},
{\bf x}, {\bf x},{\bf x},{\bf x}) + O(|| {\bf x}
||^{10}){\nonumber},
\end{eqnarray}
} \noindent where $A = f_{\bf x}(0,{\bf \mu_0})$ and {\small
\begin{equation}
B_i ({\bf x},{\bf y}) = \sum_{j,k=1}^n \frac{\partial ^2
F_i(\xi)}{\partial \xi_j \: \partial \xi_k} \bigg|_{\xi=0} x_j \;
y_k, \label{Bap}
\end{equation}
} {\small
\begin{equation}
C_i ({\bf x},{\bf y},{\bf z}) = \sum_{j,k,l=1}^n \frac{\partial ^3
F_i(\xi)}{\partial \xi_j \: \partial \xi_k \: \partial \xi_l}
\bigg|_{\xi=0} x_j \; y_k \: z_l, \label{Cap}
\end{equation}
} {\small
\begin{equation}
D_i ({\bf x},{\bf y},{\bf z},{\bf u}) = \sum_{j,k,l,r=1}^n
\frac{\partial ^4 F_i(\xi)}{\partial \xi_j \: \partial \xi_k \:
\partial \xi_l \: \partial \xi_r} \bigg|_{\xi=0} x_j \;
y_k \: z_l \: u_r, \label{Dap}
\end{equation}
} {\small
\begin{equation}
E_i ({\bf x},{\bf y},{\bf z},{\bf u},{\bf v}) = \sum_{j,k,l,r,p
=1}^n \frac{\partial ^5 F_i(\xi)}{\partial \xi_j \: \partial \xi_k
\: \partial \xi_l \: \partial \xi_r \: \partial \xi_p}
\bigg|_{\xi=0} x_j \; y_k \: z_l \: u_r \: v_p, \label{Eap}
\end{equation}
} {\small
\begin{equation}
K_i ({\bf x},{\bf y},{\bf z},{\bf u},{\bf v},{\bf w}) =
\sum_{j,\ldots,q =1}^n \frac{\partial ^6 F_i(\xi)}{\partial \xi_j \:
\partial \xi_k \: \partial \xi_l \: \partial \xi_r \: \partial \xi_p
\: \partial \xi_q} \bigg|_{\xi=0} x_j \; y_k \: z_l \: u_r \: v_p \:
w_q, \label{Kap}
\end{equation}
} {\small
\begin{equation}
L_i ({\bf x},{\bf y},{\bf z},{\bf u},{\bf v},{\bf w},{\bf t}) =
\sum_{j,\ldots,h =1}^n \frac{\partial ^7 F_i(\xi)}{\partial \xi_j
\partial \xi_k
\partial \xi_l  \partial \xi_r  \partial \xi_p  \partial \xi_q \partial \xi_h} \bigg|_{\xi=0} x_j \;
y_k \: z_l \: u_r \: v_p \: w_q \: t_h, \label{Lap}
\end{equation}
} {\small
\begin{equation}
M_i ({\bf x},{\bf y},{\bf z},{\bf u},{\bf v},{\bf w},{\bf t},{\bf
r}) = \sum_{j,\ldots,a =1}^n \frac{\partial ^8 F_i(\xi)}{\partial
\xi_j
\ldots
\partial \xi_h \partial \xi_a} \bigg|_{\xi=0}
x_j \: y_k \: z_l \: u_r \: v_p \: w_q \: t_h \: r_a, \label{Map}
\end{equation}
} {\small
\begin{equation}
N_i ({\bf x},{\bf y},{\bf z},{\bf u},{\bf v},{\bf w},{\bf t},{\bf
r},{\bf s}) = \sum_{j,\ldots,b =1}^n \frac{\partial ^9
F_i(\xi)}{\partial \xi_j
\ldots
\partial \xi_b} \bigg|_{\xi=0} x_j \;
y_k \: z_l \: u_r \: v_p \: w_q \: t_h \: r_a \: s_b, \label{Nap}
\end{equation}
} \noindent for $i = 1, \ldots, n$.

Suppose $({\bf x_0}, {\bf \mu_0})$ is an equilibrium point of
(\ref{diffequat}) where the Jacobian matrix $A$ has a pair of purely
imaginary eigenvalues $\lambda_{2,3} = \pm i \omega_0$, $\omega_0 >
0$, and admits no other eigenvalue with zero real part. Let $T^c$ be
the generalized eigenspace of $A$ corresponding to $\lambda_{2,3}$.
By this is meant that it is the largest subspace invariant by $A$ on
which the eigenvalues are $\lambda_{2,3}$.

Let $p, q \in \mathbb C ^n$ be vectors such that
\begin{equation}
A q = i \omega_0 \: q,\:\: A^{\top} p = -i \omega_0 \: p, \:\:
\langle p,q \rangle = \sum_{i=1}^n \bar{p}_i \: q_i \:\: = 1,
\label{normalization}
\end{equation}
where $A^{\top}$ is the transposed of the matrix $A$. Any vector $y
\in T^c$ can be represented as $y = w q + \bar w \bar q$, where $w =
\langle p , y \rangle \in \mathbb C$. The two dimensional center
manifold can be parameterized by $w , \bar w$, by means of an
immersion of the form  ${\bf x} = H (w, \bar w)$, where $H:\mathbb
C^2 \to \R^n$ has a Taylor expansion of the form
\begin{equation}
H(w,{\bar w}) = w q + {\bar w}{\bar q} + \sum_{2 \leq j+k \leq 9}
\frac{1}{j!k!} \: h_{jk}w^j{\bar w}^k + O(|w|^{10}), \label{defH}
\end{equation}
with $h_{jk} \in \mathbb C ^n$ and  $h_{jk}={\bar h}_{kj}$.
Substituting this expression into (\ref{diffequat}) we obtain the
following differential equation
\begin{equation} \label{ku}
H_w w' + H_{\bar w} {\bar w}' = F (H(w,{\bar w})),
\end{equation}
where $F$ is given by (\ref{Fhomo}).

The complex vectors $h_{ij}$ are obtained solving the system of
linear equations defined by the coefficients of (\ref{ku}), taking
into account the coefficients of $F$, so that system (\ref{ku}), on
the chart $w$ for a central manifold, writes as follows {\small
\[
w'= i \omega_0 w + \frac{1}{2} \; G_{21} w |w|^2 + \frac{1}{12} \;
G_{32} w |w|^4 + \frac{1}{144} \; G_{43} w |w|^6 + \frac{1}{2880} \;
G_{54} w |w|^8 + O(|w|^{10}),
\]
} with $G_{jk} \in \mathbb C$.

Solving for the vectors $h_{ij}$ the  system of linear equations
defined by the coefficients of the quadratic terms of (\ref{ku}),
taking into account the coefficients of $F$ in the expressions
(\ref{taylorexp}) and (\ref{Bap}), one has
\begin{equation}
h_{11}=-A^{-1}B(q,{\bar q}) \label{h11},
\end{equation}
\begin{equation}
h_{20}=(2i\omega_0 I_n - A)^{-1}B(q,q),\label{h20}
\end{equation}
where $I_n$ is the unit $n \times n$ matrix. Pursuing the
calculation to cubic terms, from the coefficients of the terms $w^3$
in (\ref{ku}) follows that
\begin{equation}
h_{30}=(3 i \omega_0 I_n -A)^{-1} \left[ 3B(q,h_{20})+C(q,q,q)
\right] \label{h30}.
\end{equation}

From the coefficients of the terms $w^2 {\bar w}$ in (\ref{ku}) one
obtains a singular system for $h_{21}$
\begin{equation}
(i \omega_0 I_n -A)h_{21}=C(q,q,{\bar q})+B({\bar q},h_{20})+ 2
B(q,h_{11})-G_{21}q, \label{h21m}
\end{equation}
which has a solution if and only if
\[
\langle p, C(q,q,\bar q) + B(\bar q, h_{20}) + 2 B(q,h_{11}) -G_{21}
q \rangle = 0.
\]

The {\it first Lyapunov coefficient} $l_1$ is defined by
\begin{equation}
l_1 =  \frac{1}{2} \: {\rm Re} \; G_{21}, \label{defcoef}
\end{equation}
where
\[
G_{21}= \langle p, \mathcal H_{21} \rangle, \; {\mbox {and}} \;
\mathcal H_{21} = C(q,q,\bar q) + B(\bar q, h_{20}) + 2 B(q,
h_{11}).
\]

The complex vector $h_{21}$ can be found by solving the nonsingular
$(n+1)$-dimensional system
\begin{equation}
\left( \begin{array}{cc}
i \omega_0 I_n -A & q \\
\\ {\bar p} & 0
\end{array} \right) \left( \begin{array} {c}
h_{21}\\
\\s \end{array} \right)= \left( \begin{array}{c} \mathcal H_{21} -G_{21} q \\
\\0 \end{array} \right),\label{h20}
\end{equation}
with the condition $\langle p, h_{21} \rangle = 0$. The procedure
above can be adapted in connection with the determination of
$h_{32}$ and $h_{43}$.

For the sake of completeness, in Remark \ref{nonsingular} we prove
that the system (\ref{h20}) is nonsingular and that if $(v,s)$ is a
solution of (\ref{h20}) with the condition $\langle p, v \rangle =
0$ then $v$ is a solution of (\ref{h21m}). See Remark 3.1 of
\cite{smb2}.

\begin{remark}
Write $\R^n = T^c \oplus T^{su}$, where $T^c$ and $T^{su}$ are
invariant by $A$. It can be proved that $y \in T^{su}$ if and only
if $\langle p , y \rangle = 0$. Define
\[
a = C(q,q,\bar q) + B(\bar q, h_{20}) + 2 B(q,h_{11}) -G_{21} q.
\]
Let $(v,s)$ be a solution of the homogeneous equation obtained from
(\ref{h20}). Equivalently
\begin{equation}
(i \omega_0 I_n -A)v + s q = 0, \: \: \langle p , v \rangle = 0.
\label{sistema1}
\end{equation}
From the second equation of (\ref{sistema1}), it follows that $v \in
T^{su}$, and thus $(i \omega_0 I_n -A)v \in T^{su}$. Therefore
$\langle p, (i \omega_0 I_n -A)v \rangle = 0$. Taking the inner
product of $p$ with the first equation of (\ref{sistema1}) one has
$\langle p, (i \omega_0 I_n -A)v + s q \rangle = 0$, which can be
written as $ \langle p, (i \omega_0 I_n -A)v \rangle + s \langle p,
q \rangle = 0$. Since $\langle p, q \rangle = 1$ and $\langle p, (i
\omega_0 I_n -A)v \rangle = 0$ it follows that $s = 0$. Substituting
$s = 0 $ into the first equation of (\ref{sistema1}) one has $(i
\omega_0 I_n -A)v = 0 $. This implies that
\begin{equation}
v = \alpha q, \: \alpha \in \mathbb C. \label{alpha}
\end{equation}
But $ 0 = \langle p , v \rangle = \langle p , \alpha q \rangle =
\alpha \langle p, q \rangle = \alpha$. Substituting $\alpha =0$ into
(\ref{alpha}) it follows that $v=0$. Therefore $(v,s) = (0,0)$.

Let $(v,s)$ be a solution of (\ref{h20}). Equivalently
\begin{equation}
(i \omega_0 I_n -A)v + s q = a, \: \langle p , v \rangle = 0.
\label{sistema}
\end{equation}
From the second equation of (\ref{sistema}), it follows that $v \in
T^{su}$ and thus $(i \omega_0 I_n -A)v \in T^{su}$. Therefore $
\langle p , (i \omega_0 I_n -A)v \rangle = 0$. Taking the inner
product of $p$ with the first equation of (\ref{sistema}) one has
$\langle p, (i \omega_0 I_n -A)v + s q \rangle = \langle p, a
\rangle$, which can be written as
\[
\langle p, (i \omega_0 I_n -A)v \rangle + s \langle p, q \rangle =
\langle p, a \rangle.
\]
As $\langle p, a \rangle = 0$, $\langle p, q \rangle = 1$ and
$\langle p, (i \omega_0 I_n -A)v \rangle = 0$ it follows that $s =
0$. Substituting $s = 0$ into the first equation of (\ref{sistema})
results $(i \omega_0 I_n -A)v = a$. Therefore $v$ is a solution of
(\ref{h21m}).

\label{nonsingular}
\end{remark}

From the coefficients of the terms $w^4$, $w^3 {\bar w}$ and $w^2
{\bar w ^2}$ in (\ref{ku}), one has respectively
\begin{equation}
h_{40}= (4i \omega_0 I_n -A)^{-1}[3B(h_{20},h_{20}) + 4 B(q,h_{30})+
6 C(q,q,h_{20})+ D(q,q,q,q)],\label{h40}
\end{equation}
\begin{eqnarray}\label{h31}
h_{31}=(2i \omega_0 I_n -A)^{-1}[3 B(q,h_{21}) + B(\bar q, h_{30}) +
3 B(h_{20},h_{11})  \nonumber \\ + 3 C(q,q,h_{11}) +3 C(q,\bar
q,h_{20}) + D(q,q,q,\bar q) - 3G_{21}h_{20}],
\end{eqnarray}
\begin{eqnarray}\label{h22}
h_{22}= -A^{-1} [D(q,q,\bar q, \bar q) + 4 C(q, \bar q, h_{11}) +
C(\bar q, \bar q, h_{20}) + C(q,q,{\bar h}_{20}) \nonumber \\
+ 2 B(h_{11},h_{11}) + 2 B(q,{\bar h}_{21}) + 2 B(\bar q, h_{21})+
B({\bar h}_{20},h_{20})],
\end{eqnarray}
where the term $-2h_{11}(G_{21}+{\bar G}_{21})$ has been omitted in
the last equation, since $G_{21}+{\bar G}_{21} = 0$ as $l_1 = 0$.

Defining $\mathcal H_{32}$ as
\begin{eqnarray}\label{H32}
\mathcal H_{32} = 6 B(h_{11},h_{21})+ B({\bar h}_{20},h_{30}) + 3
B({\bar h}_{21},h_{20})+ 3 B(q,h_{22}) \nonumber \\ + 2 B(\bar q,
h_{31}) +6 C(q,h_{11},h_{11}) + 3 C(q, {\bar h}_{20}, h_{20})+ 3
C(q,q,{\bar h}_{21}) \nonumber \\ +6 C(q,\bar q, h_{21}) + 6 C(\bar
q, h_{20}, h_{11}) + C(\bar q, \bar q, h_{30}) + D(q,q,q,{\bar
h}_{20}) \nonumber \\ + 6 D(q,q,\bar q,h_{11}) + 3
D(q, \bar q,\bar q, h_{20}) + E(q,q,q,\bar q,\bar q) \nonumber \\
-6 G_{21}h_{21} - 3 {\bar G}_{21} h_{21}, \nonumber
\end{eqnarray}
and from the coefficients of the terms $w^3 {\bar w}^2$ in
(\ref{ku}), one has a singular system for $h_{32}$
\begin{equation}
(i \omega_0 I_n -A)h_{32}= \mathcal H_{32} - G_{32}q, \label{h32m}
\end{equation}
which has solution if and only if
\begin{equation}
\langle p, \mathcal H_{32} - G_{32}q \rangle = 0. \label{H32m}
\end{equation}
where the terms $-6 G_{21}h_{21} - 3 {\bar G}_{21} h_{21}$ in the
last line of (\ref{H32}) actually does not enter in last equation,
since $\langle p, h_{21} \rangle = 0$.

The {\it second Lyapunov coefficient} is defined by
\begin{equation}
l_2= \frac{1}{12} \: {\rm Re} \: G_{32}, \label{defcoef2}
\end{equation}
where, from (\ref{H32m}), $G_{32}=\langle p, \mathcal H_{32}
\rangle$.

The complex vector $h_{32}$ can be found solving the nonsingular
$(n+1)$-dimensional system
\begin{equation}
\left( \begin{array}{cc}
i \omega_0 I_n -A & q \\
\\ {\bar p} & 0
\end{array} \right) \left( \begin{array} {c}
h_{32}\\
\\s \end{array} \right)= \left( \begin{array}{c} \mathcal H_{32} -G_{32} q \\
\\0 \end{array} \right),\label{h32}
\end{equation}
with the condition $\langle p, h_{32} \rangle = 0$.

From the coefficients of the terms $w^4 {\bar w}$, $w^4 {\bar w}^2$
and $w^3 {\bar w}^3$ in (\ref{ku}), one has respectively
\begin{eqnarray}\label{h41}
h_{41}=(3 i \omega_0 I_n -A)^{-1}[4B(h_{11},h_{30})+ 6
B(h_{20},h_{21}) + 4 B(q,h_{31}) \nonumber \\ + B(\bar q,h_{40}) +
12 C(q,h_{11},h_{20}) + 6 C(q,q,h_{21}) + 4 C(q,\bar q, h_{30}) \\
+ 3 C(\bar q, h_{20},h_{20}) + 4 D(q,q,q,h_{11}) + 6 D(q,q,\bar q,
h_{20}) \nonumber \\ + E(q,q,q,q,\bar q) - 6 G_{21}h_{30}]
\nonumber,
\end{eqnarray}
\begin{eqnarray}\label{h42}
h_{42}=(2 i \omega_0 I_n -A)^{-1}[8 B(h_{11},h_{31}) + 6
B(h_{20},h_{22}) + B({\bar h}_{20},h_{40}) \nonumber \\ + 6
B(h_{21},h_{21})
+ 4 B({\bar h}_{21},h_{30}) + 4 B(q,h_{32}) + 2 B(\bar q, h_{41}) \nonumber \\
+ 12 C(h_{11},h_{11},h_{20}) + 3 C(h_{20},h_{20},{\bar h}_{20}) + 24
C(q,h_{11},h_{21}) \nonumber \\ + 12 C(q,h_{20},{\bar h}_{21})
+ 4 C(q,{\bar h}_{20},h_{30}) + 6 C(q,q,h_{22}) + 8 C(q,\bar q,h_{31}) \nonumber \\
+ 8 C(\bar q, h_{11},h_{30}) + 12 C(\bar q, h_{20},h_{21}) + C(\bar
q, \bar q, h_{40}) \\ + 12 D(q,q,h_{11},h_{11})
+ 6 D(q,q,h_{20},{\bar h}_{20}) + 4 D(q,q,q,{\bar h}_{21}) \nonumber \\
+ 12 D(q,q,\bar q,h_{21}) + 24 D(q,\bar q,h_{11},h_{20}) + 4
D(q,\bar q,\bar q,h_{30}) \nonumber \\ + 3 D(\bar q,\bar q,
h_{20},h_{20}) + E(q,q,q,q,{\bar h}_{20}) + 8 E(q,q,q,\bar q,
h_{11}) \nonumber \\ + 6 E(q,q,\bar q,\bar q,h_{20}) +
K(q,q,q,q,\bar q,\bar q) \nonumber \\ - 4 (G_{32}h_{20} + 3
G_{21}h_{31} + {\bar G}_{21}h_{31})] \nonumber,
\end{eqnarray}
\begin{eqnarray}\label{h33}
h_{33}= - A^{-1}[9 B(h_{11},h_{22}) + 3 B(h_{20},{\bar h}_{31}) +
3 B({\bar h}_{20},h_{31}) + 9 B(h_{21},{\bar h}_{21}) \nonumber \\
+ B({\bar h}_{30},h_{30}) + 3 B(q,{\bar h}_{32}) + 3 B(\bar
q,h_{32}) + 6 C(h_{11},h_{11},h_{11}) \nonumber \\ + 9
C(h_{11},{\bar h}_{20},h_{20}) + 18 C(q,h_{11},{\bar h}_{21}) + 3
C(q,h_{20},{\bar h}_{30}) \nonumber \\ + 9 C(q,{\bar h}_{20},h_{21})
+ 3 C(q,q,{\bar h}_{31}) + 9 C(q,\bar q,h_{22}) + 18 C(\bar
q,h_{11},h_{21}) \nonumber \\ + 9 C(\bar q,h_{20},{\bar h}_{21}) + 3
C(\bar q, {\bar h}_{20},h_{30}) + 3 C(\bar q,\bar q, h_{31}) + 9
D(q,q,{\bar h}_{20},h_{11}) \\ + D(q,q,q,{\bar h}_{30}) + 9
D(q,q,\bar q,{\bar h}_{21}) + 18 D(q,\bar q, h_{11},h_{11})
\nonumber \\ + 9 D(q,\bar q, {\bar h}_{20}, h_{20}) + 9 D(q,\bar
q,\bar q, h_{21}) + 9 D(\bar q,\bar q, h_{11}, h_{20}) \nonumber \\
+ 3 E(q,q,q,\bar q,{\bar h}_{20}) + 9
E(q,q,\bar q,\bar q, h_{11}) + 3 E(q,\bar q,\bar q,\bar q, h_{20}) \nonumber \\
+ K(q,q,q,\bar q,\bar q, \bar q) - 3 (G_{32}+{\bar G}_{32}) h_{11} -
9 (G_{21}+{\bar G}_{21}) h_{22}] \nonumber.
\end{eqnarray}

Defining $\mathcal H_{43} $ as {\small
\begin{eqnarray}\label{H43m}
\mathcal H_{43} = 12 B(h_{11},h_{32}) +
6 B(h_{20},{\bar h}_{32}) + 3 B({\bar h}_{20},h_{41}) \nonumber \\
+ 18 B(h_{21},h_{22}) + 12 B({\bar h}_{21},h_{31}) + 4
B(h_{30},{\bar h}_{31}) + B({\bar h}_{30},h_{40}) \nonumber \\ + 4
B(q,h_{33}) + 3 B(\bar q, h_{42}) + 36 C(h_{11},h_{11},h_{21}) + 36
C(h_{11},h_{20},{\bar h}_{21}) \nonumber \\ + 12 C(h_{11},{\bar
h}_{20},h_{30}) + 3 C(h_{20},h_{20},{\bar h}_{30}) + 18
C(h_{20},{\bar h}_{20},h_{21}) \nonumber \\ + 36 C(q,h_{11},h_{22})
+ 12 C(q,h_{20},{\bar h}_{31}) + 12 C(q,{\bar h}_{20},h_{31})
\nonumber \\ + 36 C(q,h_{21},{\bar h}_{21}) + 4 C(q,h_{30},{\bar
h}_{30}) + 6 C(q,q,{\bar h}_{32}) \nonumber \\ + 12 C(q,\bar
q,h_{32}) + 24 C(\bar q,h_{11},h_{31}) + 18 C(\bar q,h_{20},h_{22})
\nonumber \\ + 3 C(\bar q,{\bar h}_{20},h_{40}) + 18 C(\bar
q,h_{21},h_{21}) + 12 C(\bar q,{\bar h}_{21},h_{30}) \nonumber \\ +
3 C(\bar q, \bar q, h_{41}) + 24 D(q,h_{11},h_{11},h_{11}) + 36
D(q,h_{11},h_{20},{\bar h}_{20}) \nonumber \\ + 36
D(q,q,h_{11},{\bar h}_{21}) + 6 D(q,q,h_{20},{\bar h}_{30}) + 18
D(q,q,{\bar h}_{20},h_{21}) \nonumber \\ + 4 D(q,q,q,{\bar h}_{31})
+ 18 D(q,q,\bar q,h_{22}) + 72 D(q,\bar q, h_{11},h_{21}) \nonumber
\\ + 36 D(q,\bar q,h_{20},{\bar h}_{21}) + 12 D(q,\bar q,{\bar
h}_{20},h_{30}) + 12 D(q,\bar q,\bar q, h_{31}) \nonumber \\ + 36
D(\bar q,h_{11},h_{11},h_{20}) + 9 D(\bar q,h_{20},h_{20},{\bar
h}_{20}) + 12 D(\bar q,\bar q,h_{11},h_{30}) \nonumber \\ + 18
D(\bar q,\bar q,h_{20},h_{21}) + D(\bar q,\bar q, \bar q,h_{40}) +
12 E(q,q,q,h_{11},{\bar h}_{20}) \nonumber \\ + E(q,q,q,q,{\bar
h}_{30}) + 12 E(q,q,q,\bar q,{\bar h}_{21}) + 36 E(q,q,\bar
q,h_{11},h_{11}) \nonumber \\ + 18 E(q,q,\bar q,h_{20},{\bar
h}_{20}) + 18 E(q,q,\bar q,\bar q, h_{21}) + 36 E(q,\bar q,\bar
q,h_{11},h_{20}) \nonumber \\ + 4 E(q,\bar q,\bar q,\bar q, h_{30})
+ 3 E(\bar q, \bar q, \bar q, h_{20},h_{20}) + 3 K(q,q,q,q,\bar
q,{\bar h}_{20}) \nonumber \\ + 12 K(q,q,q,\bar q,\bar q, h_{11}) +
6 K(q,q,\bar q,\bar q,\bar q,h_{20}) + L(q,q,q,q,\bar q,\bar q,\bar
q) \nonumber \\ - 6 (2G_{32}h_{21} + {\bar G}_{32}h_{21} + 3 G_{21}
h_{32} + 2 {\bar G}_{21}h_{32}) \nonumber,
\end{eqnarray}
} and from the coefficients of the terms $w^4 {\bar w}^3$ in
(\ref{ku}), one has a singular system for $h_{43}$
\begin{eqnarray}\label{h43m}
(i \omega_0 I_n -A)h_{43}= \mathcal H_{43}- G_{43} q
\end{eqnarray}
which has solution if and only if
\begin{eqnarray}\label{h43}
\langle p, \mathcal H_{43}- G_{43} q \rangle =0.
\end{eqnarray}
where the terms $- 6 (2G_{32}h_{21} + {\bar G}_{32}h_{21} + 3 G_{21}
h_{32} + 2 {\bar G}_{21}h_{32})$ appearing in the last line of
equation (\ref{H43m}) actually do not enter in the last equation,
since $\langle p, h_{21} \rangle = 0$ and $\langle p, h_{32} \rangle
= 0$.

The {\it third Lyapunov coefficient} is defined by
\begin{equation}
l_3= \frac{1}{144} \: {\rm Re} \: G_{43}, \label{defcoef3}
\end{equation}
where, from (\ref{h43}), $ G_{43} = \langle p, \mathcal H_{43}
\rangle$.

The complex vector $h_{43}$ can be found solving the nonsingular
$(n+1)$-dimensional system
\begin{equation}
\left( \begin{array}{cc}
i \omega_0 I_n -A & q \\
\\ {\bar p} & 0
\end{array} \right) \left( \begin{array} {c}
h_{43}\\
\\s \end{array} \right)= \left( \begin{array}{c} \mathcal H_{43} -G_{43} q \\
\\0 \end{array} \right),\label{h32}
\end{equation}
with the condition $\langle p, h_{43} \rangle = 0$.

Defining $\mathcal H_{54} $ by the expression below
{\footnotesize
\begin{eqnarray*}
20 B(h_{11},h_{43}) + 10 B(h_{20},\bar h_{43}) + 6 B(\bar
h_{20},h_{52}) + 40 B(h_{21},h_{33}) +30 B(\bar h_{21},h_{42}) + \\
60 B(h_{22},h_{32}) +  10 B(h_{30},\bar h_{42}) + 4 B(\bar
h_{30},h_{51}) + 40 B(h_{31},\bar h_{32}) + 20 B(\bar h_{31},h_{41})
+ \\ 5 B(h_{40},\bar h_{41}) + B(\bar h_{40},h_{50}) + 5 B(q,h_{44})
+ 4 B(\bar q,h_{53}) + 120 C(h_{11},h_{11},h_{32}) + \\ 60
C(h_{11},\bar h_{20},h_{41}) + 360 C(h_{11},h_{21},h_{22}) + 240
C(h_{11},\bar h_{21},h_{31}) +  80 C(h_{11},h_{30},\bar h_{31}) + \\
20 C(h_{11},\bar h_{30},h_{40}) + 120 C(h_{20},h_{11},\bar h_{32}) +
15 C(h_{20},h_{20},\bar h_{41}) + 60 C(h_{20},\bar h_{20},h_{32}) +
\\ 120 C(h_{20},h_{21},\bar h_{31}) + 180 C(h_{20},\bar
h_{21},h_{22}) + 10 C(h_{20},h_{30},\bar h_{40}) + 40 C(h_{20},\bar
h_{30},h_{31}) + \\ 3 C(\bar h_{20},\bar h_{20},h_{50}) + 120 C(\bar
h_{20},h_{21},h_{31}) + 30 C(\bar h_{20},\bar h_{21},h_{40}) + 60
C(\bar h_{20},h_{30},h_{22}) + \\ 180 C(h_{21},h_{21},\bar h_{21}) +
60 C(\bar h_{21},\bar h_{21},h_{30}) + 40 C(h_{30},h_{21},\bar
h_{30}) + 80 C(q,h_{11},h_{33}) + \\ 30 C(q,h_{20},\bar h_{42}) + 30
C(q,\bar h_{20},h_{42}) + 120 C(q,h_{21},\bar h_{32}) + 120 C(q,\bar
h_{21},h_{32}) + \\ 90 C(q,h_{22},h_{22}) + 20 C(q,h_{30},\bar
h_{41}) + 20 C(q,\bar h_{30},h_{41}) + 80 C(q,h_{31},\bar h_{31}) +
5 C(q,h_{40},\bar h_{40}) + \\ 10 C(q,q,\bar h_{43}) + 20 C(q,\bar
q,h_{43}) + 60 C(\bar q,h_{11},h_{42}) + 40 C(\bar q,h_{20},h_{33})
+ 12 C(\bar q,\bar h_{20},h_{51}) + \\ 120 C(\bar q,h_{21},h_{32}) +
60 C(\bar q,\bar h_{21},h_{41}) + 40 C(\bar q,h_{30},\bar h_{32}) +
4 C(\bar q,\bar h_{30},h_{50}) + \\ 120 C(\bar q,h_{31},h_{22}) + 20
C(\bar q,h_{40},\bar h_{31}) + 6 C(\bar q,\bar q,h_{52}) + 240
D(h_{11},h_{11},h_{11},h_{21}) + \\ 120 D(h_{11},h_{11},\bar
h_{20},h_{30}) + 360 D(h_{20},h_{11},h_{11},\bar h_{21}) + 360
D(h_{20},h_{11},\bar h_{20},h_{21}) + \\ 60
D(h_{20},h_{20},h_{11},\bar h_{30}) + 90 D(h_{20},h_{20},\bar
h_{20},\bar h_{21}) + 30 D(h_{20},\bar h_{20},\bar h_{20},h_{30}) +
\\ 360 D(q,h_{11},h_{11},h_{22}) + 240 D(q,h_{11},\bar
h_{20},h_{31}) + 720 D(q,h_{11},h_{21},\bar h_{21}) + \\ 80
D(q,h_{11},h_{30},\bar h_{30}) + 240 D(q,h_{20},h_{11},\bar h_{31})
+ 15 D(q,h_{20},h_{20},\bar h_{40}) + \\ 180 D(q,h_{20},\bar
h_{20},h_{22}) + 120 D(q,h_{20},h_{21},\bar h_{30}) + 180
D(q,h_{20},\bar h_{21},\bar h_{21}) + \\ 15 D(q,\bar h_{20},\bar
h_{20},h_{40}) + 180 D(q,\bar h_{20},\bar h_{20},h_{21}) + 120
D(q,\bar h_{20},h_{30},\bar h_{21}) + \\ 120 D(q,q,h_{11},\bar
h_{32}) + 30 D(q,q,h_{20},\bar h_{41}) + 60 D(q,q,\bar
h_{20},h_{32}) + 120 D(q,q,h_{21},\bar h_{31}) + \\ 180 D(q,q,\bar
h_{21},h_{22}) + 10 D(q,q,h_{30},\bar h_{40}) + 40 D(q,q,\bar
h_{30},h_{31}) + 10 D(q,q,q,\bar h_{42}) + \\ 40 D(q,q,\bar
q,h_{33}) + 240 D(q,\bar q,h_{11},h_{32}) + 120 D(q,\bar
q,h_{20},\bar h_{32}) + 60 D(q,\bar q,\bar h_{20},h_{41}) + \\ 360
D(q,\bar q,h_{21},h_{22}) + 240 D(q,\bar q,\bar h_{21},h_{31}) + 80
D(q,\bar q,h_{30},\bar h_{31}) + 20 D(q,\bar q,\bar h_{30},h_{40}) +
\\ 30 D(q,\bar q,\bar q,h_{42}) + 240 D(\bar q,h_{11},h_{11},h_{31})
+ 60 D(\bar q,h_{11},\bar h_{20},h_{40}) + 360 D(\bar
q,h_{11},h_{21},h_{21}) + \\ 240 D(\bar q,h_{11},h_{30},\bar h_{21})
+ 360 D(\bar q,h_{20},h_{11},h_{22}) + 60 D(\bar
q,h_{20},h_{20},\bar h_{31}) +
\\ 120 D(\bar q,h_{20},\bar h_{20},h_{31}) + 360 D(\bar
q,h_{20},h_{21},\bar h_{21}) + 40 D(\bar q,h_{20},h_{30},\bar
h_{30}) +  \\ 120 D(\bar q,\bar h_{20},h_{30},h_{21}) + 60 D(\bar
q,\bar q,h_{11},h_{41}) + 60 D(\bar q,\bar q,h_{20},h_{32}) + 6
D(\bar q,\bar q,h_{20},h_{50}) + \\ 120 D(\bar q,\bar
q,h_{21},h_{31}) + 30 D(\bar q,\bar q,\bar h_{1},h_{40}) + 60 D(\bar
q,\bar q,h_{30},h_{22}) + 4 D(\bar q,\bar q,\bar q,h_{51}) + \\ 120
E(q,h_{11},h_{11},h_{11},h_{11}) + 360 E(q,h_{20},h_{11},h_{11},\bar
h_{20}) + 45 E(q,h_{20},h_{20},\bar h_{20},\bar h_{20}) + \\ 360
E(q,q,h_{11},h_{11},\bar h_{21}) + 360 E(q,q,h_{11},\bar
h_{20},h_{21}) + 120 E(q,q,h_{20},h_{11},\bar h_{30}) + \\ 180
E(q,q,h_{20},\bar h_{20},\bar h_{21}) + 30 E(q,q,\bar h_{20},\bar
h_{20},h_{30}) + 80 E(q,q,q,h_{11},\bar h_{31}) + \\ 10
E(q,q,q,h_{20},\bar h_{40}) + 60 E(q,q,q,\bar h_{20},h_{22}) + 40
E(q,q,q,h_{21},\bar h_{30}) + 60 E(q,q,q,\bar h_{21},\bar h_{21}) +
\\ 5 E(q,q,q,q,\bar h_{41}) + 40 E(q,q,q,\bar q,\bar h_{32}) + 360
E(q,q,\bar q,h_{11},h_{22}) + 120 E(q,q,\bar q,h_{20},\bar h_{31}) +
\\ 120 E(q,q,\bar q,\bar h_{20},h_{31}) + 360 E(q,q,\bar
q,h_{21},\bar h_{21}) + 40 E(q,q,\bar q,h_{30},\bar h_{30}) + 60
E(q,q,\bar q,\bar q,h_{32}) + \\ 720 E(q,\bar
q,h_{11},h_{11},h_{21}) + 240 E(q,\bar q,h_{11},\bar h_{20},h_{30})
+ 720 E(q,\bar q,h_{20},h_{11},\bar h_{21}) + \\ 60 E(q,\bar
q,h_{20},h_{20},\bar h_{30}) + 360 E(q,\bar q,h_{20},\bar
h_{20},h_{21}) + 240 E(q,\bar q,\bar q,h_{11},h_{31}) + \\ 180
E(q,\bar q,\bar q,h_{20},h_{22}) + 30 E(q,\bar q,\bar q,\bar
h_{20},h_{40}) + 180 E(q,\bar q,\bar
q,h_{21},h_{21}) + 120 E(q,\bar q,\bar q,h_{30},\bar h_{21}) + \\
20 E(q,\bar q,\bar q,\bar q,h_{41}) + 240 E(\bar
q,h_{20},h_{11},h_{11},h_{11}) + 180 E(\bar
q,h_{20},h_{20},h_{11},\bar h_{20}) + \\ 120 E(\bar q,\bar
q,h_{11},h_{11},h_{30}) + 360 E(\bar q,\bar q,h_{20},h_{11},h_{21})
+ 90 E(\bar q,\bar q,h_{20},h_{20},\bar h_{21}) + \\ 60 E(\bar
q,\bar q,h_{20},\bar h_{20},h_{30}) + 20 E(\bar q,\bar q,\bar
q,h_{11},h_{40}) + 40 E(\bar q,\bar q,\bar
q,h_{20},h_{31}) + 40 E(\bar q,\bar q,\bar q,h_{30},h_{21}) + \\
E(\bar q,\bar q,\bar q,\bar q,h_{50}) + 120
K(q,q,q,h_{11},h_{11},\bar h_{20}) + 30 K(q,q,q,h_{20},\bar
h_{20},\bar h_{20}) + \\ 20 K(q,q,q,q,h_{11},\bar h_{30}) + 30
K(q,q,q,q,\bar h_{20},\bar h_{21}) + K(q,q,q,q,q,\bar h_{40}) + 20
K(q,q,q,q,\bar q,\bar h_{31}) + \\ 240 K(q,q,q,\bar q,h_{11},\bar
h_{21}) + 40 K(q,q,q,\bar q,h_{20},\bar h_{30}) + 120 K(q,q,q,\bar
q,\bar h_{20},h_{21}) + \\ 60 K(q,q,q,\bar q,\bar q,h_{22}) + 240
K(q,q,\bar q,h_{11},h_{11},h_{11}) + 360 K(q,q,\bar
q,h_{20},h_{11},\bar h_{20}) + \\ 360 K(q,q,\bar q,\bar
q,h_{11},h_{21}) + 180 K(q,q,\bar q,\bar q,h_{20},\bar h_{21}) + 60
K(q,q,\bar q,\bar q,\bar h_{20},h_{30}) + \\ 40 K(q,q,\bar q,\bar
q,\bar q,h_{31}) + 360 K(q,\bar q,\bar q,h_{20},h_{11},h_{11}) + 90
K(q,\bar q,\bar q,h_{20},h_{20},\bar h_{20}) + \\ 80 K(q,\bar q,\bar
q,\bar q,h_{11},h_{30}) + 120 K(q,\bar q,\bar q,\bar
q,h_{20},h_{21}) + 5 K(q,\bar q,\bar q,\bar q,\bar q,h_{40}) + \\ 60
K(\bar q,\bar q,\bar q,h_{20},h_{20},h_{11}) + 10 K(\bar q,\bar
q,\bar q,\bar q,h_{20},h_{30}) + 3 L(q,q,q,q,q,\bar h_{20},\bar
h_{20}) + \\ 4 L(q,q,q,q,q,\bar q,\bar h_{30}) + 60 L(q,q,q,q,\bar
q,h_{11},\bar h_{20}) + 30 L(q,q,q,q,\bar q,\bar q,\bar h_{21}) + \\
120 L(q,q,q,\bar q,\bar q,h_{11},h_{11}) + 60 L(q,q,q,\bar q,\bar
q,h_{20},\bar h_{20}) + 40 L(q,q,q,\bar q,\bar q,\bar q,h_{21}) +
\\ 120 L(q,q,\bar q,\bar q,\bar q,h_{20},h_{11}) + 10 L(q,q,\bar q,\bar q,\bar q,\bar
q,h_{30}) + 15 L(q,\bar q,\bar q,\bar q,\bar q,h_{20},h_{20}) + \\
6 M(q,q,q,q,q,\bar q,\bar q,\bar h_{20}) + 20 M(q,q,q,q,\bar q,\bar
q,\bar q,h_{11}) + 10 M(q,q,q,\bar q,\bar q,\bar q,\bar q,h_{20}) +
\\ N(q,q,q,q,q,\bar q,\bar q,\bar q,\bar q),
\end{eqnarray*}
} and from the coefficients of the terms $w^5 {\bar w}^4$ in
(\ref{ku}), one has a singular system for $h_{54}$
\begin{eqnarray}\label{h54m}
(i \omega_0 I_n -A)h_{54}= \mathcal H_{54}- G_{54} q
\end{eqnarray}
which has solution if and only if
\begin{eqnarray}\label{h54}
\langle p, \mathcal H_{54}- G_{54} q \rangle =0.
\end{eqnarray}

The {\it fourth Lyapunov coefficient} is defined by
\begin{equation}
l_4= \frac{1}{2880} \: {\rm Re} \: G_{54}, \label{defcoef4}
\end{equation}
where, from (\ref{h54}), $ G_{54} = \langle p, \mathcal H_{54}
\rangle$.

\begin{remark}\label{conceitual}
Other equivalent definitions and algorithmic  procedures to write
the expressions for the Lyapunov coefficients $l_j , j= 1,2,3,4$,
for two dimensional systems can be found in Andronov et al.
\cite{al} and Gasull et al. \cite{gt}, among others. These
procedures apply also to the $n$--dimensional systems of this work,
if properly restricted to the center manifold. The authors found,
however, that the  method  outlined above, due to Kuznetsov
\cite{kuznet, kuznet2}, requiring  no explicit formal evaluation of
the center manifold, is better adapted to the needs of long
calculations in \cite{smb1, smb2, smb3}, and for that in
\cite{smb4}, where $n=3$.
\end{remark}

A {\it Hopf point} $({\bf x_0}, {\bf \mu_0})$ is an equilibrium
point of (\ref{diffequat}) where the Jacobian matrix $A = f_{\bf
x}({\bf x_0}, {\bf \mu_0})$ has a pair of purely imaginary
eigenvalues $\lambda_{2,3} = \pm i \omega_0$, $\omega_0 > 0$, and
admits  no other critical eigenvalues ---i.e. located on the
imaginary axis. At a Hopf point a two dimensional center manifold is
well-defined, it is invariant under the flow generated by
(\ref{diffequat}) and can be continued with arbitrary high class of
differentiability to nearby parameter values. In fact, what is well
defined is the $\infty$-jet ---or infinite Taylor series--- of the
center manifold, as well as that of its continuation, any two of
them having contact in the arbitrary high  order of their
differentiability class.

A Hopf point is called {\it transversal} if the parameter dependent
complex eigenvalues cross the imaginary axis with non-zero
derivative. In a neighborhood of a transversal Hopf point
---H1 point, for concision--- with $l_1 \neq 0$ the dynamic
behavior of the system (\ref{diffequat}), reduced to the family of
parameter-dependent continuations of the center manifold, is
orbitally topologically equivalent to the following complex normal
form
\[
w' = (\eta + i \omega) w + l_1 w |w|^2 ,
\]
$w \in \mathbb C $, $\eta$, $\omega$ and $l_1$ are real functions
having  derivatives of arbitrary  high order, which are
continuations  of $0$, $\omega_0$ and the first Lyapunov coefficient
at the H1 point. See  \cite{kuznet}. As $l_1 < 0$ ($l_1 > 0$) one
family of stable (unstable) periodic orbits can be found on this
family of manifolds, shrinking  to an equilibrium point at the H1
point.

A {\it Hopf point of codimension 2} is a Hopf point where $l_1$
vanishes. It is called {\it transversal} if $\eta = 0$ and $l_1 = 0$
have transversal intersections, where $\eta = \eta (\mu)$ is the
real part of the critical eigenvalues. In a neighborhood of a
transversal Hopf point of codimension 2 ---H2 point, for
concision---  with $l_2 \neq 0$ the dynamic behavior of the system
(\ref{diffequat}), reduced to the family of parameter-dependent
continuations of the center manifold, is orbitally topologically
equivalent to
\[
w' = (\eta + i \omega_0) w + \tau w |w|^2 + l_2 w |w|^4 ,
\]
where $\eta$ and $\tau$ are unfolding parameters.  See
\cite{kuznet}. The bifurcation diagrams for $l_2 \neq 0$ can be
found in \cite{kuznet}, p. 313, and in \cite{takens}.

A {\it Hopf point of codimension 3} is a Hopf point of codimension 2
where $l_2$ vanishes. A Hopf point of codimension 3 point is called
{\it transversal} if $\eta = 0$, $l_1 = 0$ and $l_2 = 0$ have
transversal intersections. In a neighborhood of a transversal Hopf
point of codimension 3 ---H3 point, for concision--- with $l_3 \neq
0$ the dynamic behavior of the system (\ref{diffequat}), reduced to
the family of parameter-dependent continuations of the center
manifold, is orbitally topologically equivalent to
\[
w' = (\eta + i \omega_0) w + \tau w |w|^2 + \nu w |w|^4 + l_3 w
|w|^6 ,
\]
where $\eta$, $\tau$ and $\nu$ are unfolding parameters. The
bifurcation diagram for $l_3 \neq 0$ can be found in Takens
\cite{takens} and in \cite{smb3}.

A {\it Hopf point of codimension 4} is a Hopf point of codimension 3
where $l_3$ vanishes. A Hopf point of codimension 4 is called {\it
transversal} if $\eta = 0$, $l_1 = 0$, $l_2 = 0$ and $l_3 = 0$ have
transversal intersections. In a neighborhood of a transversal Hopf
point of codimension 4 ---H4 point, for concision--- with $l_4 \neq
0$ the dynamic behavior of the system (\ref{diffequat}), reduced to
the family of parameter-dependent continuations of the center
manifold, is orbitally topologically equivalent to
\[
w' = (\eta + i \omega_0) w + \tau w |w|^2 + \nu w |w|^4 + \sigma w
|w|^6 + l_4 w |w|^8 ,
\]
where $\eta$, $\tau$, $\nu$ and $\sigma$ are unfolding parameters.

\begin{teo}
Suppose that the system
\[
{\bf x}' = f({\bf x},{\bf \mu}), \: {\bf x}=(x,y,z), \: \mu =
(\beta, \alpha, \kappa, \varepsilon)
\]
has the equilibrium ${\bf x} = {\bf 0}$ for $\mu = 0$ with
eigenvalues
\[
\lambda_{2,3} (\mu) = \eta (\mu) \pm i \omega(\mu),
\]
where $\omega(0) = \omega_0 > 0$. For $\mu = 0$ the following
conditions hold
\[
\eta (0) = 0, \: l_1(0) = 0, \: l_2(0) = 0, \: l_3(0) = 0,
\]
where $l_1(\mu)$, $l_2(\mu)$ and $l_3(\mu)$ are the first, second
and third Lyapunov coefficients, respectively. Assume that the
following genericity conditions are satisfied
\begin{enumerate}
\item $l_4 (0) \neq 0$, where $l_4 (0)$ is the fourth Lyapunov
coefficient;

\item the map $\mu \to (\eta (\mu), l_1 (\mu), l_2 (\mu), l_3
(\mu))$ is regular at $\mu = 0$.

\end{enumerate}
Then, by the introduction of a complex variable, the above system
reduced to the family of parameter-dependent continuations of the
center manifold, is orbitally topologically equivalent to
\[
w' = (\eta + i \omega_0) w + \tau w |w|^2 + \nu w |w|^4 + \sigma w
|w|^6 + l_4 w |w|^8
\]
where $\eta$, $\tau$, $\nu$ and $\sigma$ are unfolding parameters.

\label{teoremaHopf}
\end{teo}

\begin{remark}\label{needs}
The expressions for the Lyapunov coefficients in this article will
be of interest in the study of the Hopf bifurcations of codimensions
1, 2, 3 and 4 in the Watt governor system with a spring \cite{smb4},
pursuing previous bifurcation analysis carried out in \cite{smb3}.

\end{remark}

\vspace{0.2cm}

\noindent {\bf Acknowledgement}: The first and second authors
developed this work under the projects CNPq grants 473824/04-3 and
473747/2006-5. The first author is fellow of CNPq. The third author
is supported by CAPES. This work was finished while the second
author visited Universitat Aut\`onoma de Barcelona, supported by
CNPq grant 210056/2006-1.


\begin{thebibliography}{99}

\bibitem{al} A. A. Andronov,  E. A. Leontovich et al.,
{\it Theory of Bifurcations of Dynamic Systems on a Plane}, Halsted
Press, J. Wiley \& Sons, New York, 1973.

\bibitem{gt} A. Gasull and J. Torregrosa, A new approach to the computation
of the Lyapunov Constants, Comp. and Appl. Math., {\bf 20} (2001),
149--177.

\bibitem{kuznet} Y. A. Kuznetsov, {\it Elements of Applied
Bifurcation Theory}, Springer--Verlag, New York, 2004.

\bibitem{kuznet2} Y. A. Kuznetsov, Numerical normalization
techniques for all codim 2 bifurcations of equilibria in ODE's, SIAM
J. Numer. Anal., {\bf 36} (1999), 1104--1124.

\bibitem{smb1} J. Sotomayor, L. F. Mello and D. C. Braga, Stability
and Hopf bifurcation in the Watt governor system, Commun. Appl.
Nonlinear Anal. {\bf 13} (2006), 4, 1--17.

\bibitem{smb2} J. Sotomayor, L. F. Mello and D. C. Braga, Bifurcation analysis of the
Watt governor system, Comp. Appl. Math. {\bf 26} (2007), 19--44.

\bibitem{smb3} J. Sotomayor, L. F. Mello and D. C. Braga, Stability and Hopf bifurcation
in an hexagonal governor system, Nonlinear Anal.: Real World Appl.
(2007), in press,  doi: 10.1016/nonrwa.2007.01.007.

\bibitem{smb4} J. Sotomayor, L. F. Mello and D. C. Braga, Hopf bifurcations in a Watt
governor with a spring, Preprint (2007).

\bibitem{takens} F. Takens, Unfoldings of certain singularities of
vectorfields: Generalized Hopf bifurcations, J. Diff. Equat., {\bf
14} (1973), 476--493.

\bibitem{mello} Site with the files used in computer assited
arguments in this work:
  {\small http://www.ici.unifei.edu.br/luisfernando/wgss}

\bibitem{math} S. Wolfram, {\it The Mathematica Book}, fifth
edition, Wolfram Media Inc., Champaign (2003).

\end{thebibliography}
\end{document}